\documentclass[11pt,a4paper]{article}
\usepackage[T1]{fontenc}
\usepackage[utf8]{inputenc}
\usepackage[margin=2.6cm]{geometry}
\usepackage{amsmath,amssymb,amsthm,mathtools}
\usepackage{enumitem,booktabs,array}
\usepackage{xcolor}
\usepackage[colorlinks=true,linkcolor=blue!60!black,urlcolor=blue!60!black,citecolor=blue!60!black]{hyperref}

\newcommand{\Q}{\mathbb{Q}}
\newcommand{\C}{\mathbb{C}}
\newcommand{\Vol}{\operatorname{Vol}}
\newcommand{\La}{\Lambda}
\newcommand{\lhi}[1]{\log_{\ge n}(#1)}
\newcommand{\sgn}{\operatorname{sgn}}
\newcommand{\pr}[2]{(#1,#2)}
\newcommand{\Dehn}{\mathcal{D}}
\newcommand{\Lg}{\mathcal{L}}
\newcommand{\disc}{\operatorname{disc}}
\newcommand{\Gram}{\operatorname{Gram}}

\theoremstyle{plain}
\newtheorem{theorem}{Theorem}[section]
\newtheorem{proposition}[theorem]{Proposition}
\newtheorem{lemma}[theorem]{Lemma}

\newtheorem*{theoremA}{Theorem A}
\theoremstyle{definition}
\newtheorem{definition}[theorem]{Definition}
\newtheorem{rem}[theorem]{Remark}

\title{The infinitesimal volume of a tetrahedron}
\author{Sinan \"Unver}
\date{}

\begin{document}
\maketitle

\begingroup
\NoHyper
\renewcommand{\thefootnote}{}
\endNoHyper
\endgroup

\begin{abstract}
Let $k_n$ denote the truncated polynomial ring of length $n$ over a field $k$ of characteristic $0$. We construct an infinitesimal volume for tetrahedra over $k_n$ with respect to a non-degenerate quadratic form (e.g.\ spherical or hyperbolic tetrahedra). It takes values in the first Andr\'e--Quillen homology $D_1(k_n)$ of $k_n$ relative to $\Q$, up to a twist by the orientation, and it satisfies the scissors congruence relations. The construction involves lifting the tetrahedron to a ring of twice the length. The main point is that the volume does not depend on this lift, which we deduce from the classical Schl\"afli formula for spherical tetrahedra. The infinitesimal volume is modelled on the additive dilogarithm on the Bloch group.
\end{abstract}

\section*{Introduction}

The scissors congruence group $\mathcal P(S^3)$ of spherical tetrahedra has two basic invariants: the volume, and the Dehn invariant $\sum_\ell\theta_\ell\otimes\psi_\ell$. Here $\theta_\ell$ is the length of the edge $\ell$ and $\psi_\ell$ is the dihedral angle along it. By the work of Dupont, Parry and Sah, the hyperbolic and spherical Dehn complexes in dimension $3$ are closely related to the Bloch complex
\[
  B_2(\C)\xrightarrow{\ \delta\ }\Lambda^2\C^\times,\qquad \delta[x] := (1-x)\wedge x,
\]
whose kernel is $K_3^{\rm ind}(\C)$ and whose cokernel is $K_2(\C)$, up to torsion \cite{DS}, \cite{DPS}, \cite[\S6.1]{G1}. In the hyperbolic case, Dupont and Sah relate the Dehn invariant to the exact sequence of Bloch and Wigner \cite[Proposition 4.5, Theorem 4.10]{DS}, and show that the volume of the ideal tetrahedron with vertices $\infty,0,1,z$ is the Bloch--Wigner dilogarithm $D(z)$ \cite[(4.13)]{DS}.

Goncharov gave a motivic interpretation of these facts. A non-Euclidean geodesic simplex defines a framed mixed Hodge--Tate structure; the Dehn invariant becomes its coproduct \cite[Theorem 1.5]{G1}, and the volume is one of its periods \cite[(5) and Theorem 4.7]{G1}. In this language the Schl\"afli formula is the differential equation satisfied by the periods \cite[\S4.8]{G1}. Goncharov defined a Hopf algebra $S_\bullet(\C)$, generated by simplices relative to a non-degenerate quadric \cite[\S3.4]{G1}. Its definition makes sense over any field $F$, and for $F$ algebraically closed he conjectured a homomorphism $s_2:S_2(F)\to B_2(F)$ compatible with the Dehn invariant and $\delta$ \cite[Conjecture 6.2]{G1}.

Over a field $k$ of characteristic $0$ there is no $k$-valued volume or dilogarithm. Over the truncated polynomial rings $k_n := k[t]/(t^n)$, however, there are infinitesimal versions of both. Cathelineau's additive Bloch group $\beta_2(k)$ is defined by the four-term functional equation of the entropy function. It carries the infinitesimal Dehn invariant $\langle a\rangle\mapsto a\otimes a+(1-a)\otimes(1-a)$, whose cokernel is $\Omega^1_k$ \cite[Th\'eor\`eme 1]{Ca}; see also \cite{Ca2}. Goncharov proposed to view Euclidean geometry as a degeneration of hyperbolic geometry over the dual numbers \cite[\S6.3]{G1}. The Euclidean scissors congruence groups make sense over any field \cite[\S3.2]{G2}. They carry a volume map, and the analogue of Sydler's theorem predicts that it identifies the kernel of the Euclidean Dehn invariant with the base field, over $\mathbb R$ \cite[Conjecture 6.3]{G1} and over algebraically closed fields \cite[Conjecture 1.1]{G2}.

Bloch and Esnault constructed an additive dilogarithm using the localization sequence in $K$-theory \cite{BE}. In \cite{U1}, the additive dilogarithm
\[
  Li_{2,n}:B_2(k_n)\to k^{\oplus(n-1)}
\]
was constructed on the Bloch group itself. Its components have $t$-degrees $n+1,\dots,2n-1$, and it induces an isomorphism $\ker(\delta)^\circ\simeq K_3(k_n,(t))^{(2)}_\Q\simeq k^{\oplus(n-1)}$ (Theorem 1.3.2 and Corollary 5.1.3 there). As observed in \S1.3 there, this is an analogue of Sydler's theorem: $\delta$ plays the role of the Dehn invariant and the components of $Li_{2,n}$ play the role of volumes. Volume maps of this kind are also expected in higher weight \cite[\S2.2]{U3}.

In this paper we construct an infinitesimal volume for tetrahedra over $k_n$ with respect to a non-degenerate quadratic form. Let $V$ be a free $k_n$-module of rank four with a non-degenerate symmetric bilinear form $Q$. An \emph{orientation} of $V$ is a linear form $\omega$ on $\Lambda^4V$ whose square is the Gram determinant: $\omega(v_1\wedge\dots\wedge v_4)^2 = \det\big(Q(v_j,v_{j'})\big)$. The constant term of the Gram determinant of a basis of $V$ lies in $k^\times$, and its class $\disc(V)\in k^\times/k^{\times2}$ does not depend on the basis; we write $\sqrt{\disc(V)}$ for a square root of a representative. Orientations exist over $k'_n$, where $k' := k\big(\sqrt{\disc(V)}\big)$, and they are defined over $k_n$ if and only if $\disc(V)$ is trivial. A \emph{tetrahedron} in $\mathbb P(V)$ is a quadruple of points of $\mathbb P(V)$ over $k_n$, any one, two or three of which span a non-degenerate subspace (\S\ref{ssec:quad}). The basic case is $V = k_n^4$ with the standard form and the orientation $\det$, which gives spherical tetrahedra; the form $x_1^2+x_2^2+x_3^2-x_4^2$, with the orientation $i\det$, gives hyperbolic ones.

In place of the Dehn invariant we use a \emph{Dehn element}, built from one elementary invariant. Let $W$ be a plane with a non-degenerate quadratic form $Q$ and an orientation, i.e.\ a linear form $\omega$ on $\Lambda^2W$ with $\omega(v\wedge w)^2 = Q(v,v)Q(w,w)-Q(v,w)^2$. Then $Q(Jv,w) = \omega(v\wedge w)$ defines a rotation $J$ by a right angle, $J^2 = -1$, and after adjoining $i = \sqrt{-1}$ the plane splits into the eigenlines $W^\pm$ of $J$ for $\pm i$. These are the two isotropic lines; their images $q^\pm$ in $\mathbb P(W)$ are the two isotropic points, and the orientation is what tells them apart. For non-isotropic points $x_j = [v_j]$ of $\mathbb P(W)$, with $v_j = v_j^++v_j^-$, put
\[
  \langle x_0,x_1\rangle_{W,\omega} := \frac{Q(v_0^-,v_1^+)}{Q(v_0^+,v_1^-)} = \frac{Q(v_0,v_1)+i\,\omega(v_0\wedge v_1)}{Q(v_0,v_1)-i\,\omega(v_0\wedge v_1)} .
\]
This is the cross-ratio $r(x_0,x_1,q^+,q^-)$ of $x_0$, $x_1$ and the two isotropic points. Laguerre's formula says that the angle from $x_0$ to $x_1$ is $\frac1{2i}\log r(x_0,x_1,q^+,q^-)$; indeed, for two unit vectors in the Euclidean plane at angle $\theta$ we get $e^{2i\theta}$. For an edge $\ell$ of a tetrahedron we apply this twice: to the two vertices on the line $\ell$, which gives $\mathcal Z_\ell$, and to the two faces through $\ell$, regarded as points of the pencil of planes through $\ell$, which gives $\mathcal W_\ell$. For a real spherical tetrahedron, $\mathcal Z_\ell = \exp(2i\theta_\ell)$ and $\mathcal W_\ell = \exp(\pm2i\psi_\ell)$, where $\theta_\ell$ is the length of $\ell$ and $\psi_\ell$ the dihedral angle along it. The Dehn element is
\[
  \Dehn(\Delta) := \sum_\ell\sgn(\ell)\,\mathcal Z_\ell\wedge\mathcal W_\ell\ \in\ \Lambda^2E_m^\times,
\]
where $\sgn(\ell) = \pm1$ is the sign of the permutation listing the vertices of $\ell$ and then the other two, each in increasing order,
$E_m := E[t]/(t^m)$, and $E$ is a finite extension of $k$ containing $\sqrt{-1}$ and the square roots needed to orient the planes of the edges (Definition~\ref{def:dehn}). For $a,b\in k_{2n}^\times$ let
\[
  \La_n(a\wedge b) := d\log(a)\otimes\log_{\ge n}(b)-d\log(b)\otimes\log_{\ge n}(a)\ \in\ \Omega^1_{k_n/k}\otimes_{k_n}t^nk_{2n},
\]
where $\log_{\ge n}$ is the part of $t$-degree $\ge n$ of $\log^\circ(a) := \log(a/a(0))$. The target of $\La_n$ is isomorphic to $t^n\Omega^1_{k_{2n}/k}$, and $d$ identifies the first Andr\'e--Quillen homology $D_1(k_n) = t^{n+1}k_{2n}$ of $k_n$ relative to $\Q$ with it.

\begin{theoremA}
Let $k$ be a field of characteristic $0$, $n\ge2$, $(V,Q)$ a quadratic space of rank four over $k_n$, and $\omega$ an orientation of $V\otimes_{k_n}k'_n$. For a tetrahedron $\Delta$ in $\mathbb P(V)$, choose $E\supseteq k'$ as above, and a lift of $(V,Q,\omega,\Delta)\otimes_{k_n}E_n$ to $E_{2n}$: an oriented quadratic space $(\tilde V,\tilde Q,\tilde\omega)$ over $E_{2n}$ with a tetrahedron $\tilde\Delta$ in $\mathbb P(\tilde V)$, whose reduction modulo $t^n$ is identified with $(V,Q,\omega,\Delta)\otimes_{k_n}E_n$. Put
\[
  \Vol^\omega_n(\Delta) := -\frac14\,d^{-1}\circ\La_n\big(\Dehn(\tilde\Delta)\big)\ \in\ D_1(E_n).
\]
\begin{enumerate}[nosep]
  \item[(i)] $\Vol^\omega_n(\Delta)$ does not depend on the choice of the lift, of $E$, or of the orientations of the edges, and it lies in $\sqrt{\disc(V)}\cdot D_1(k_n)$. In particular it lies in $D_1(k_n)$ if $\disc(V)$ is trivial.
  \item[(ii)] $\Vol^{-\omega}_n = -\Vol^\omega_n$, and $\Vol^{\omega'}_n(g\Delta) = \Vol^\omega_n(\Delta)$ for an isometry $g:V\to V'$ onto another quadratic space which carries an orientation $\omega'$ of $V'$ to $\omega$. Moreover $\Vol^\omega_n$ satisfies the scissors congruence relations (R1)--(R4) of \S\ref{ssec:vol}, and so induces a homomorphism
  \[
    \Vol^\omega_n:\ \mathcal P(V)\longrightarrow D_1(k'_n)
  \]
  on the scissors congruence group of tetrahedra in $\mathbb P(V)$.
  \item[(iii)] $\Vol^\omega_n$ is natural for the automorphisms $\sigma_\lambda:t\mapsto\lambda t$, $\lambda\in k^\times$, of $k_n$: the infinitesimal volume of the base change of $(V,Q,\omega,\Delta)$ along $\sigma_\lambda$ is $\sigma_\lambda\big(\Vol^\omega_n(\Delta)\big)$, where $\sigma_\lambda$ acts on $D_1(k_n)$ by $\lambda^w$ on $t^w$. In particular, for $V = k_n^4$ with the standard form, the component of $t$-degree $w$, $n+1\le w\le2n-1$, is homogeneous of degree $w$ in the jets of the vertices; for $n = 2$ the only $t$-degree is $3$, the degree of homogeneity of the Euclidean volume.
\end{enumerate}
For $V = k_n^4$ with the standard form and $\omega = \det$ this gives the infinitesimal volume $\Vol_n:\mathcal P(k_n)\to D_1(k_n)$ of spherical tetrahedra; for $x_1^2+x_2^2+x_3^2-x_4^2$ and $\omega = i\det$ it gives an infinitesimal volume of hyperbolic tetrahedra, with values in $i\cdot D_1(k_n)$.
\end{theoremA}

By (ii), $\Vol^\omega_n(\Delta)\otimes\omega$ does not depend on the orientation, and by (i) it is invariant under $\operatorname{Gal}(k'/k)$. So it is an element of $D_1(k_n)\otimes_{k_n}(\Lambda^4V)^\vee$, defined for every quadratic space of rank four, oriented or not: this is what we mean by a twist by the orientation. Similarly, the volume of a Euclidean simplex with vertices $0,v_1,v_2,v_3$ is naturally the element $\frac16v_1\wedge v_2\wedge v_3$ of the determinant line, and becomes a number only once an orientation is chosen.

The theorem is proved in Section~2 (Theorems~\ref{thm:vol}, \ref{thm:scissors} and~\ref{thm:descent}). We call $\Vol^\omega_n(\Delta)$ the \emph{infinitesimal volume} of $\Delta$. The factor $-\frac14$ makes it, formally, the integral of the truncated Schl\"afli form $\sum_\ell(d\theta_\ell\otimes(\psi_\ell)_{\ge n}-d\psi_\ell\otimes(\theta_\ell)_{\ge n})$. It is not identically zero; see the examples in Remarks~\ref{rem:n2} and~\ref{rem:jets}.

The proof of (i) has two steps. First, the variation of $\La_n(f\wedge g)$ under a change of lift of a $k_n$-point is the pairing of $d\log f\wedge d\log g$ with a derivation (Lemma~\ref{lem:var}). Second, the relevant $2$-form vanishes. For the Dehn element this $2$-form is $\sum_\ell\sgn(\ell)\,d\log\mathcal Z_\ell\wedge d\log\mathcal W_\ell$, which for real spherical tetrahedra is $-4\sum_\ell d\theta_\ell\wedge d\psi_\ell = -4\,d\big(\sum_\ell\theta_\ell\,d\psi_\ell\big)$, the derivative of Schl\"afli's formula $2\,dV = \sum_\ell\theta_\ell\,d\psi_\ell$. The orientations of the edges drop out of this $2$-form, so it is a rational $2$-form on the space of configurations, and its vanishing follows from the classical case (Proposition~\ref{prop:schlafli}). The independence of $E$ and of the orientations of the edges follows from a Galois argument (Theorem~\ref{thm:descent}).

The model for this construction is the additive dilogarithm. If we replace the Dehn element by $\delta[x] = (1-x)\wedge x$, the $2$-form becomes $d\log(1-x)\wedge d\log x$, which vanishes on $\mathbb A^1\setminus\{0,1\}$, and $d^{-1}\circ\La_n\circ\delta$ gives a renormalization of the additive dilogarithm of \cite{U1} on $B_2(k_n)$ (Remark~\ref{rem:dilog}). Up to sign, $\La_n\circ\delta$ is a truncated, infinitesimal version of the $1$-form $\log(1-z)\,d\log z-\log z\,d\log(1-z)$, which Dupont and Sah use, following Bloch, to relate the Bloch--Wigner dilogarithm to the volume of ideal hyperbolic tetrahedra \cite[(4.14)--(4.18)]{DS}.

The paper is organized as follows. After fixing notation, we define $\La_n$ in Section~1 and compute how it changes under a change of lift. In Section~2 we introduce the angle between two points of an oriented quadratic plane, use it to define the Dehn element of a tetrahedron, prove the Schl\"afli identity, and construct the infinitesimal volume over $E$; in \S\ref{ssec:descent} we descend it to $k$ and prove Theorem~A. The final remarks discuss explicit formulas, examples, the hyperbolic case, and the expected map from tetrahedra over $k_n$ to $B_2(k_n)$, the analogue of Goncharov's conjectural homomorphism $s_2$.

\section*{Notation and conventions}

Let $k$ be a field of characteristic $0$. For $m\ge1$ let $k_m := k[t]/(t^m)$. An element of $k_m$ has \emph{$t$-degree} $w$ if it lies in $kt^w$. The group $k^\times$ acts on $k_m$ by the $k$-algebra automorphisms $t\mapsto\lambda t$, and $kt^w$ is the subspace on which $\lambda$ acts by $\lambda^w$. For $a\in k_m$ we write $a(0)\in k$ for its constant term.

For a functor $F$ on $k$-algebras of the form $k_m$, we let $F^\circ := F(k_m)^\circ := \ker\big(F(k_m)\to F(k)\big)$ denote its \emph{infinitesimal part}. Since $k\to k_m\to k$ is the identity, $F(k_m) = F(k)\oplus F(k_m)^\circ$. For instance $k_m^\circ = tk_m$ and $(k_m^\times)^\circ = 1+tk_m$. For a unit $a$ of $k_m$ we put $\log^{\circ}(a) := \log(a/a(0))\in k_m^\circ$. The map $\log^\circ:k_m^\times\to k_m^\circ$ is a homomorphism: it is the projection $a\mapsto a/a(0)$ onto $(k_m^\times)^\circ$, followed by the logarithm $(k_m^\times)^\circ\simeq k_m^\circ$.

For a commutative ring $A$ we let $A^\flat$ denote the set of units $a\in A^\times$ such that $1-a$ is also a unit. The \emph{Bloch group} $B_2(A)$ is the $\Q$-vector space generated by the symbols $[a]$, $a\in A^\flat$, subject to the five-term relations
\[
  [x]-[y]+[y/x]-\Big[\frac{1-x^{-1}}{1-y^{-1}}\Big]+\Big[\frac{1-x}{1-y}\Big] = 0,
\]
for $x,y\in A^\times$ with $(1-x)(1-y)(1-x/y)\in A^\times$. We put $\Lambda^2A^\times := (\Lambda^2_{\mathbb Z}A^\times)\otimes\Q$. The map $[a]\mapsto\delta[a] := (1-a)\wedge a$ induces the differential $\delta: B_2(A)\to\Lambda^2A^\times$ of the Bloch complex. If $A$ is local, then $A^\flat$, and hence the five-term relations, only depend on the residue field. For instance, $a\in k_m^\flat$ if and only if $a(0)\in k\setminus\{0,1\}$. We write $\ker(\delta)^\circ := \ker(\delta)\cap B_2(k_m)^\circ$.

Let $R$ be a local ring and $L$ a free $R$-module of rank $2$. The \emph{cross-ratio} of four points $x_j = [v_j]$ of $\mathbb P(L)$, $j = a,b,c,d$, is
\[
  r(x_a,x_b,x_c,x_d) := \frac{\mu(v_a\wedge v_c)\,\mu(v_b\wedge v_d)}{\mu(v_a\wedge v_d)\,\mu(v_b\wedge v_c)}\ \in\ R,
\]
where $\mu:\Lambda^2L\to R$ is any isomorphism; it is defined when the denominator is a unit, and does not depend on $\mu$ or on the representatives $v_j$. In an affine coordinate $z$ it is $\frac{(z_a-z_c)(z_b-z_d)}{(z_a-z_d)(z_b-z_c)}$. Equivalently, if $f$ is a coordinate on $\mathbb P(L)$ with a pole at $x_c$ and a zero at $x_d$, then
\begin{equation}\label{eq:cr-f}
  r(x_a,x_b,x_c,x_d) = \frac{f(x_b)}{f(x_a)} .
\end{equation}
In particular $r(x_a,x_a,x_c,x_d) = 1$ and $r(x_a,x_b,x_c,x_d)\,r(x_b,x_e,x_c,x_d) = r(x_a,x_e,x_c,x_d)$.

As in \cite[\S3.1]{U2}, Andr\'e--Quillen homology is taken relative to $\Q$: for a commutative $\Q$-algebra $A$, $D_1(A) := D_1(A|\Q)$ denotes the first Andr\'e--Quillen homology of $A$ over $\Q$ \cite[\S3.5]{Lo}. Let $A = P/J$ with $P$ a polynomial algebra over $k$. In the transitivity sequence
\[
  D_1(P|\Q,A)\to D_1(A)\to D_1(A|P)\to\Omega^1_{P/\Q}\otimes_PA
\]
for $\Q\to P\to A$, the first term vanishes since $P$ is a filtered colimit of smooth $\Q$-algebras, and $D_1(A|P) = J/J^2$, the last map being induced by $d$. Hence
\begin{equation}\label{eq:D1}
  D_1(A) = \ker\big(J/J^2\xrightarrow{\ d\ }\Omega^1_{P/\Q}\otimes_PA\big).
\end{equation}

\section{\texorpdfstring{The map $\Lambda_n$}{The map Lambda n}}\label{sec:La}

For a $k$-algebra $A$ and $a \in A^{\times}$ we put $d\log (a):=da/a \in \Omega ^{1} _{A/k}$. We will use that $\Omega^1_{k_m/k} = k_m\,dt/(t^{m-1}dt)$, because $d(t^m) = mt^{m-1}dt$ and $m$ is invertible in $k$.

Let $1\le n<m$. Every $a\in k_m$ can be written uniquely as $a = a_{<n}+a_{\ge n}$, with $a_{<n}$ in the $k$-span of $1,t,\dots,t^{n-1}$ and $a_{\geq n} \in t^n k_m$. Let $\log _{\geq n}(a):= (\log ^{\circ}(a))_{\geq n}$; this is a homomorphism $k_m^\times\to t^nk_m$. Define
\[
\La_n:\ \Lambda^2k_{2n}^\times\longrightarrow\Omega^1_{k_n/k}\otimes_{k_n}t^nk_{2n},\qquad
a\wedge b\longmapsto d\log(a)\otimes\lhi{b}-d\log(b)\otimes\lhi{a}.
\]
Here $d\log(a)$ is the logarithmic differential of the image of $a$ in $k_n$. The formula is biadditive and alternating. Since the target is a $k$-vector space, $\La_n$ is defined on $\Lambda^2k_{2n}^\times = (\Lambda^2_{\mathbb Z}k_{2n}^\times)\otimes\Q$.

\subsection{Variation under a change of lift}

Let $B$ be a $k$-algebra, and let $v, v' \in {\rm Spec} (B)(k_{2n})$ be two $k_{2n}$-valued points with the same reduction $\overline{v} \in {\rm Spec} (B)(k_{n})$. Let $\varphi, \varphi': B \to k_{2n}$ and $\overline{\varphi}: B \to k_{n}$ be the corresponding $k$-algebra maps. The ideal $t^n k_{2n}$ is annihilated by $t^n$, so it is a $k_{n}$-module, and hence a $B$-module via $\overline{\varphi}$. With this module structure, $D:=\varphi'-\varphi : B \to t^nk_{2n}$ is a $k$-derivation. Indeed, $D$ is $k$-linear, and for $b,c\in B$,
\[
  D(bc) = \big(\varphi'(b)-\varphi(b)\big)\varphi'(c)+\varphi(b)\big(\varphi'(c)-\varphi(c)\big) = \overline\varphi(c)\,D(b)+\overline\varphi(b)\,D(c).
\]
Let $\iota_{D}: \Omega ^{1} _{B/k} \to t^nk_{2n}$ be the $B$-linear map corresponding to $D$, and extend it to $\iota _{D}: \Omega ^{2} _{B/k} \to  \Omega^{1} _{B/k} \otimes _{B} t^nk_{2n}$ by
\[
\iota_{D}(\eta _{1} \wedge \eta_{2}):= \eta _{1} \otimes \iota _{D}(\eta _{2})  -\eta_{2} \otimes \iota _{D}(\eta _{1}) .
\]

\begin{lemma}\label{lem:var}
With the above notation, for $f, g \in B^{\times}$ we have
\begin{equation}\label{eq:var}
  \La_n(f(v')\wedge g(v'))-\La_n(f(v)\wedge g(v)) = \overline{\varphi } _{*}\big(\iota_{D}(d\log f\wedge d\log g)\big),
\end{equation}
where $\overline{\varphi }_*:\Omega^1_{B/k}\otimes_{B} t^nk_{2n} \to\Omega^1_{k_n/k}\otimes_{k_n} t^n k_{2n}$ is the map induced by $\overline {\varphi}$. In particular, if $d\log f\wedge d\log g = 0$ in $\Omega^2_{B/k}$, then $\La_n(f(v)\wedge g(v))$ only depends on $\overline v$.
\end{lemma}

\begin{proof}
Since $v$ and $v'$ have the same reduction $\overline v$, the differentials $d\log f(v)$ and $d\log f(v')$ in $\Omega^1_{k_n/k}$ are both equal to $d\log f(\overline v)$, and similarly for $g$. So only the terms $\log_{\ge n}$ change. We have $g(v')=g(v)+D(g)$ with $D(g)\in t^nk_{2n}$. In particular $g(v')$ and $g(v)$ have the same constant term, so
\[
  \log^\circ(g(v'))-\log^\circ(g(v)) = \log\big(1+D(g)/g(v)\big) = D(g)/g(v),
\]
since $(t^nk_{2n})^2 = 0$. This lies in $t^nk_{2n}$, and equals $D(g)/g(\overline v)$ for the $k_n$-module structure. Hence
\[
  d \log f(\overline {v}) \otimes \log _{\geq n} g(v')-d \log f(\overline {v}) \otimes \log _{\geq n} g(v)
  =\overline{\varphi}_{*}(d \log f \otimes \iota _{D}(d \log g)).
\]
In the same way,
\[
  d \log g(\overline {v}) \otimes \log _{\geq n} f(v')-d \log g(\overline {v}) \otimes \log _{\geq n} f(v)
  =\overline{\varphi}_{*}(d \log g \otimes \iota _{D}(d \log f)).
\]
Subtracting the second identity from the first gives \eqref{eq:var}.
\end{proof}

\subsection{The target}

Let $P = k[t]\to k_n$ be the presentation with kernel $J = (t^n)$. By \eqref{eq:D1},
\[
  D_1(k_n) = \ker\big(J/J^2\xrightarrow{d}\Omega^1_{P/\Q}\otimes_Pk_n\big)=t^{n+1}k_{2n}.
\]
Indeed, $J/J^2 = t^nk_{2n}$ and $\Omega^1_{P/\Q}\otimes_Pk_n = (\Omega^1_{k/\Q}\otimes_kk_n)\oplus k_n\,dt$. For $f = \sum_{w=n}^{2n-1}c_wt^w$ we have $d f = \sum_{w}\big(wc_wt^{w-1}dt+t^w\,dc_w\big)$; the terms $t^w\,dc_w$ vanish in $\Omega^1_{P/\Q}\otimes_Pk_n$ since $w\ge n$, and the rest vanishes in $k_n\,dt$ if and only if $c_n = 0$. So $D_1(k_n)$ has basis $t^w$, $n+1\le w\le 2n-1$. On the other hand, $t^n\Omega^1_{k_{2n}/k}$ has basis $t^{w-1}dt$, $n+1\le w\le 2n-1$, and $d(t^w) = w\,t^{w-1}dt$. Hence $d:t^{n+1}k_{2n}\to t^n\Omega^1_{k_{2n}/k}$ is an isomorphism, with inverse
\begin{equation}\label{eq:dinv}
  d^{-1}\big(t^{w-1}dt\big) = \frac{t^w}{w}\qquad(n+1\le w\le 2n-1),
\end{equation}
the antiderivative without constant term. In particular $D_1(k_1) = 0$, and from now on $n\ge2$.

The target of $\La_{n}$ is identified with $t^n \Omega ^{1} _{k_{2n}/k}$ by the map $\eta \otimes f\mapsto f \tilde{\eta}$, where $\tilde{\eta} \in \Omega ^{1} _{k_{2n}/k}$ is any lift of $\eta$. This is well defined: the kernel of $\Omega^1_{k_{2n}/k}\to\Omega^1_{k_n/k}$ is $t^{n-1}\Omega^1_{k_{2n}/k}$, and $t^n k_{2n}\cdot t^{n-1}\Omega^1_{k_{2n}/k} = 0$. It is bijective, because the source has $k$-basis $t^idt\otimes t^n$, $0\le i\le n-2$, which is mapped to the basis $t^{n+i}dt$ of $t^n\Omega^1_{k_{2n}/k}$. We obtain a map
\[
  d^{-1}\circ\La_{n}: \Lambda ^{2} k_{2n} ^{\times} \to D_{1}(k_n).
\]

\begin{lemma}\label{lem:explicit}
Let $a,b\in k_{2n}^\times$ with $\log^\circ(a) = \sum_i\alpha_it^i$ and $\log^\circ(b) = \sum_j\beta_jt^j$. Then
\[
  d^{-1}\circ\La_n(a\wedge b) = \sum_{w=n+1}^{2n-1}\frac1w\Big(\sum_{1\le i\le w-n}i\,(\alpha_i\beta_{w-i}-\beta_i\alpha_{w-i})\Big)t^w .
\]
\end{lemma}

\begin{proof}
In $\Omega^1_{k_n/k}$ we have $d\log(a) = \sum_{1\le i\le n-1}i\alpha_it^{i-1}dt$, and $\log_{\ge n}(b) = \sum_{j\ge n}\beta_jt^j$. Using $t^{2n-1}dt = 0$ in $\Omega^1_{k_{2n}/k}$, we get
\[
  d\log(a)\otimes\log_{\ge n}(b) = \sum_{\substack{1\le i\le n-1\\ j\ge n}}i\alpha_i\beta_j\,t^{i+j-1}dt
  = \sum_{w=n+1}^{2n-1}\Big(\sum_{1\le i\le w-n}i\,\alpha_i\beta_{w-i}\Big)t^{w-1}dt,
\]
since for $i+j = w\le 2n-1$ the condition $j\ge n$ is $i\le w-n$, which implies $i\le n-1$. Subtract the same expression with $a$ and $b$ interchanged and apply \eqref{eq:dinv}.
\end{proof}

\begin{rem}\label{rem:dilog}
Let $B := k[z,\frac{1}{z(1-z)}]$. Since $\Omega^2_{B/k} = 0$, Lemma~\ref{lem:var}, applied to $f := 1-z$ and $g := z$, shows that $d^{-1}\circ\La_n(\delta[x])$, for $x\in k_{2n}^\flat$, only depends on the image of $x$ in $k_n$. Since the five-term relations only involve constant terms, it descends to a map $\Lg: B_2(k_n)\to D_1(k_n)$. By Lemma~\ref{lem:explicit} and \cite[Proposition 2.1.2, Definition 2.2.4]{U1}, $\Lg = -\sum_w\frac1w Li_{2,n,w}\,t^w$, where $Li_{2,n} = \oplus_wLi_{2,n,w}$ is the additive dilogarithm of \cite{U1}, defined with the differential $[x]\mapsto x\wedge(1-x)$. Along a lift of the form $s\exp(u)$, the formula for $\Lg$ is the one in \cite[\S3]{U5}, and a map of the same shape as $\La_n$, in weight three, is used in \S6.1 there. We do not use $\Lg$ below; it is the model for the construction of the infinitesimal volume, with $\delta[x]$ replaced by the Dehn element.
\end{rem}

\section{The infinitesimal volume of a tetrahedron}\label{sec:vol}

The Dehn element is built from one elementary invariant: the angle between two points of a projective line equipped with a quadratic form and an orientation (\S\ref{ssec:angles}). Applied to an edge of a tetrahedron, it gives the length of the edge; applied to the link of the edge, i.e.\ to the two faces through it, it gives the dihedral angle (\S\ref{ssec:dehn}). The angle involves $\sqrt{-1}$, and orientations involve square roots. We therefore first work over a field $E\supseteq k$ which contains $i := \sqrt{-1}$ and the square roots that are needed (\S\S\ref{ssec:quad}--\ref{ssec:vol}), and then descend to $k$ (\S\ref{ssec:descent}).

For a field $E\supseteq k$ we put $E_m := E[t]/(t^m)$. Section~1 applies with $k$ replaced by $E$; in particular $D_1(E_n) = t^{n+1}E_{2n} = D_1(k_n)\otimes_kE$.

\subsection{Oriented quadratic spaces and simplices}\label{ssec:quad}

Let $R$ be a local ring in which $2$ is invertible; we shall only use $R = E_m$ and $R = \mathbb R$. For vectors $v_1,\dots,v_s$ of an $R$-module with a symmetric bilinear form $Q$, we write
\[
  \Gram(v_1,\dots,v_s) := \det\big(Q(v_j,v_{j'})\big)_{1\le j,j'\le s}
\]
for their \emph{Gram determinant}. Replacing the $v_j$ by unit multiples multiplies it by a square, so whether it is a unit only depends on the lines $Rv_j$. A \emph{quadratic space} over $R$ is a free $R$-module $V$ of finite rank $r$ with a symmetric bilinear form $Q$ which is non-degenerate, i.e.\ $\Gram(e_1,\dots,e_r)$ is a unit for a basis $e_1,\dots,e_r$. An \emph{orientation} of $V$ is a linear form $\omega:\Lambda^rV\to R$ with $\omega(e_1\wedge\dots\wedge e_r)^2 = \Gram(e_1,\dots,e_r)$ for one, and hence every, basis. In other words, $Q$ induces a form on the line $\Lambda^rV$, namely $(v_1\wedge\dots\wedge v_r,w_1\wedge\dots\wedge w_r)\mapsto\det\big(Q(v_j,w_{j'})\big)$, and an orientation is an isometry from this line onto $R$ with the standard form $(x,y)\mapsto xy$. If an orientation $\omega$ exists, then there are exactly two of them, $\omega$ and $-\omega$. For $R = E_m$ and $R = \mathbb R$, one exists if and only if the image of the Gram determinant in the residue field is a square. Choosing a basis, $V = R^r$ and $Q(v,w) = \pr{v}{w}_A := v^{\mathsf T}Aw$ for a symmetric matrix $A$ with $\det A\in R^\times$; then an orientation is $v_1\wedge\dots\wedge v_r\mapsto\omega\det(v_1,\dots,v_r)$ for some $\omega\in R$ with $\omega^2 = \det A$, and we also call $\omega$ an orientation of $A$.

A \emph{non-degenerate subspace} of $V$ is a submodule $W$ spanned by vectors $w_1,\dots,w_s$ with $\Gram(w_1,\dots,w_s)$ a unit. Then $w_1,\dots,w_s$ is a basis of $W$, $(W,Q|_W)$ is a quadratic space, $V = W\oplus W^\perp$, and the quotient $V/W\simeq W^\perp$ is a quadratic space of rank $r-s$, with the form induced from $W^\perp$. Given orientations $\omega$ of $V$ and $\omega_W$ of $W$, there is a unique orientation $\omega_{V/W}$ of $V/W$, the \emph{induced orientation}, such that
\begin{equation}\label{eq:induced}
  \omega(w\wedge u) = \omega_W(w)\,\omega_{V/W}(\bar u)\qquad(w\in\Lambda^sW,\ u\in\Lambda^{r-s}V),
\end{equation}
where $\bar u$ is the image of $u$ in $\Lambda^{r-s}(V/W)$. Indeed, the left hand side only depends on $\bar u$, and computing with a basis of $W$ and a basis of $W^\perp$ shows that the Gram determinant of $V$ is the product of those of $W$ and of $V/W$, so that $\omega_{V/W}$ is an orientation. Replacing $\omega_W$ by $-\omega_W$ replaces $\omega_{V/W}$ by $-\omega_{V/W}$.

The points of $\mathbb P(V)$ over $R$ are the lines $Rv$ spanned by vectors $v$ which are part of a basis; we write $[v]$. A point $[v]$ is \emph{non-isotropic} if $Q(v,v)\in R^\times$. A \emph{simplex} in $\mathbb P(V)$ is a sequence $\Delta = ([v_1],\dots,[v_r])$ of $r = \operatorname{rank}V$ points such that for every proper subset $F\subsetneq\{1,\dots,r\}$ the vectors $v_j$, $j\in F$, span a non-degenerate subspace $W_F$. (The $v_j$ themselves need not span $V$. We call $\Delta$ \emph{degenerate} if one of the $v_j$ is a linear combination of the others.) For $1\le|F|<r$ we have the \emph{face} $\Delta_F := ([v_j])_{j\in F}$, a simplex in $\mathbb P(W_F)$, and the \emph{link} $\Delta^F := ([\bar v_j])_{j\notin F}$, a sequence of points of $\mathbb P(V/W_F)$; in both, the indices are taken in increasing order. We write $\sgn(F)$ for the sign of the permutation which lists $F$ and then its complement, each in increasing order.

\begin{lemma}\label{lem:link}
Let $\Delta$ be a simplex in $\mathbb P(V)$, $F\subsetneq\{1,\dots,r\}$ nonempty, and $G,G'\subseteq\{1,\dots,r\}\setminus F$ of the same size. Let $\beta$ be the minor of the Gram matrix $\big(Q(v_j,v_{j'})\big)$ with rows $F\cup G$ and columns $F\cup G'$, ordered as $F$ followed by $G$, resp.\ by $G'$. Then
\[
  \beta = \det\big(Q(v_j,v_{j'})\big)_{j,j'\in F}\cdot\det\big(Q(\bar v_g,\bar v_{g'})\big)_{g\in G,\,g'\in G'} .
\]
In particular the link $\Delta^F$ is a simplex in $\mathbb P(V/W_F)$.
\end{lemma}

\begin{proof}
Replace each $v_g$, $g\in G\cup G'$, by its projection $v_g'$ to $W_F^\perp$, which represents $\bar v_g$. This subtracts from the last rows and columns of the minor combinations of the first $|F|$ ones, so $\beta$ is unchanged. Since $Q(v_j,v_g') = 0$ for $j\in F$, the minor becomes block triangular, and the formula follows. For $G = G'$ with $F\cup G$ proper, the left hand side is a unit, which shows that the link is a simplex.
\end{proof}

For a simplex of dimension three, i.e.\ $r = 4$, the faces of size $2$ are the edges, and the link of an edge $\ell = \{a,b\}$ consists of the images $\bar v_c$, $\bar v_d$ of the other two vertices in the plane $V/W_\ell$. Geometrically, $\mathbb P(W_\ell)$ is the line through $[v_a]$ and $[v_b]$, and $\mathbb P(V/W_\ell)$ is the pencil of planes through this line; $[\bar v_c]$ and $[\bar v_d]$ are the two faces of the tetrahedron containing the edge. This is the algebraic version of the familiar fact that the dihedral angle along an edge is the angle between the two faces in the plane orthogonal to the edge.

\subsection{The angle between two points of an oriented plane}\label{ssec:angles}

Let $(V,Q,\omega)$ be an oriented quadratic space of rank $2$ over $R$. Thus $\omega$ is a linear form on $\Lambda^2V$ subject to Lagrange's identity
\[
  \omega(v\wedge w)^2 = Q(v,v)\,Q(w,w)-Q(v,w)^2\qquad(v,w\in V).
\]
Since $Q$ is non-degenerate, there is a unique endomorphism $J$ of $V$ with
\[
  Q(Jv,w) = \omega(v\wedge w)\qquad(v,w\in V).
\]
It is the rotation by a right angle: for $V = \mathbb R^2$ with the standard form and orientation, $J(x,y) = (-y,x)$.

\begin{lemma}\label{lem:J}
$J$ is an isometry of $V$, and $J^2 = -1$.
\end{lemma}

\begin{proof}
As $\omega$ is alternating, $J$ is skew-adjoint: $Q(Jv,w) = -Q(v,Jw)$. In a basis $e_1,e_2$ with Gram matrix $G$, and with $s := \omega(e_1\wedge e_2)$, the matrix of $J$ is $-sG^{-1}\varepsilon$, where $\varepsilon = \left(\begin{smallmatrix}0&1\\-1&0\end{smallmatrix}\right)$. Its trace vanishes, $G^{-1}$ being symmetric and $\varepsilon$ alternating, and its determinant is $s^2/\det G = 1$. Hence $J^2 = -1$ by Cayley--Hamilton, and $Q(Jv,Jw) = -Q(v,J^2w) = Q(v,w)$.
\end{proof}

Assume now that $R$ contains $i$. Then $V = V^+\oplus V^-$, where $V^\pm$ is the eigenline of $J$ for the eigenvalue $\pm i$, and we write $v = v^++v^-$ accordingly; explicitly $v^\pm = \frac12(v\mp iJv)$. Both lines are isotropic: for $x,y\in V^+$ we have $Q(x,y) = Q(Jx,Jy) = i^2Q(x,y)$, so $Q(x,y) = 0$, and similarly for $V^-$. Hence
\[
  Q(v,w) = Q(v^+,w^-)+Q(v^-,w^+),
\]
$Q$ restricts to a perfect pairing $V^+\times V^-\to R$, and the images $q^\pm$ of $V^\pm$ in $\mathbb P(V)$ are the two isotropic points. Replacing $\omega$ by $-\omega$ replaces $J$ by $-J$ and interchanges $V^+$ and $V^-$; so does replacing $i$ by $-i$. In other words, once $i$ is fixed, an orientation is the same as an ordering of the two isotropic points. If $v$ is non-isotropic, then $Q(v,v) = 2Q(v^+,v^-)$ is a unit, so $v^+$ and $v^-$ are generators of $V^+$ and $V^-$.

\begin{definition}\label{def:angle}
For non-isotropic points $x_0 = [v_0]$ and $x_1 = [v_1]$ of $\mathbb P(V)$, the \emph{angle} from $x_0$ to $x_1$ is
\[
  \langle x_0,x_1\rangle_{V,\omega} := \frac{Q(v_0^-,v_1^+)}{Q(v_0^+,v_1^-)}\ \in\ R^\times.
\]
For a simplex $\Delta = (x_0,x_1)$ in $\mathbb P(V)$ we write $\langle\Delta\rangle_{V,\omega}$. When $V$ and $\omega$ are clear from the context we drop them from the notation, and we also write $\langle v_0,v_1\rangle$ for $\langle[v_0],[v_1]\rangle$.
\end{definition}

The angle does not change if $v_0$ or $v_1$ is multiplied by a unit, so it only depends on the points. To see what it is, choose generators $e_\pm$ of $V^\pm$ and write $v^\pm = \varphi_\pm(v)\,e_\pm$. For non-isotropic $v$ put $f(v) := \varphi_+(v)/\varphi_-(v)$; it is a unit, and it only depends on $[v]$. Since $Q(v_0^-,v_1^+) = \varphi_-(v_0)\varphi_+(v_1)\,Q(e_-,e_+)$, we get
\begin{equation}\label{eq:angle-f}
  \langle x_0,x_1\rangle_{V,\omega} = \frac{f(x_1)}{f(x_0)} .
\end{equation}
Thus $f$ is a coordinate on $\mathbb P(V)\simeq\mathbb P^1$ with its zero at $q^-$ and its pole at $q^+$, well defined up to a constant factor, and the angle is the ratio of its values.

\begin{lemma}\label{lem:angle}
Let $x_0,x_1,x_2$ be non-isotropic points of $\mathbb P(V)$, $x_j = [v_j]$.
\begin{enumerate}[nosep]
  \item[(a)] $\displaystyle\langle x_0,x_1\rangle_{V,\omega} = \frac{Q(v_0,v_1)+i\,\omega(v_0\wedge v_1)}{Q(v_0,v_1)-i\,\omega(v_0\wedge v_1)}$, and the product of the numerator and the denominator is $Q(v_0,v_0)\,Q(v_1,v_1)$.
  \item[(b)] $\langle x_1,x_0\rangle_{V,\omega} = \langle x_0,x_1\rangle_{V,\omega}^{-1}$, and replacing $\omega$ by $-\omega$, or $i$ by $-i$, replaces $\langle x_0,x_1\rangle_{V,\omega}$ by its inverse.
  \item[(c)] $\langle x_0,x_1\rangle_{V,\omega}\,\langle x_1,x_2\rangle_{V,\omega} = \langle x_0,x_2\rangle_{V,\omega}$.
  \item[(d)] If $g:V\to V'$ is an isometry onto an oriented quadratic space $(V',Q',\omega')$ carrying $\omega'$ to $\omega$, then $\langle gx_0,gx_1\rangle_{V',\omega'} = \langle x_0,x_1\rangle_{V,\omega}$.
\end{enumerate}
\end{lemma}

\begin{proof}
(a) We have $Q(v_0,v_1)+i\omega(v_0\wedge v_1) = Q(v_0+iJv_0,v_1) = 2Q(v_0^-,v_1) = 2Q(v_0^-,v_1^+)$, and similarly $Q(v_0,v_1)-i\omega(v_0\wedge v_1) = 2Q(v_0^+,v_1^-)$. The product is $Q(v_0,v_1)^2+\omega(v_0\wedge v_1)^2 = Q(v_0,v_0)Q(v_1,v_1)$ by Lagrange's identity. (b) and (c) follow from \eqref{eq:angle-f}, together with the fact that $\pm\omega$ and $\pm i$ interchange $V^+$ and $V^-$, which inverts $f$. (d) $g$ carries $J$, $V^\pm$ and $Q$ to the corresponding data on $V'$.
\end{proof}

\begin{rem}\label{rem:euclid}
Let $V = \mathbb R^2$ with the standard form and orientation, extended to $\C$. Then $V^\pm = \C\cdot(1,\mp i)$, and for $v = (\cos\theta,\sin\theta)$ we have $v^\pm = \frac12e^{\pm i\theta}(1,\mp i)$, so that $f(v) = e^{2i\theta}$. Hence if $v_0,v_1$ are unit vectors and $v_1$ is at angle $\theta$ from $v_0$, then $\langle v_0,v_1\rangle = e^{2i\theta}$. The factor $2$ reflects the fact that the angle only depends on points of $\mathbb P^1$, where angles are defined modulo $\pi$.
\end{rem}

\begin{rem}\label{rem:laguerre}
Since $f$ has its pole at $q^+$ and its zero at $q^-$, \eqref{eq:angle-f} and \eqref{eq:cr-f} give $\langle x_0,x_1\rangle_{V,\omega} = r(x_0,x_1,q^+,q^-)$, the cross-ratio of $x_0$, $x_1$ and the two isotropic points. This is Laguerre's formula: the angle from $x_0$ to $x_1$ is $\frac1{2i}\log r(x_0,x_1,q^+,q^-)$. In other words, $\mathbb P(V)$ minus the isotropic points is a torsor under the multiplicative group, the orientation tells us which isotropic point is $0$ and which is $\infty$, and angles are differences of logarithms. Equivalently, the angle is the integral of the invariant $1$-form $\omega(x\wedge dx)/Q(x,x) = \frac1{2i}\,d\log f$ from $x_0$ to $x_1$. Similarly, up to sign, the volume of a spherical tetrahedron $T\subseteq\mathbb P^3(\mathbb R)$ is the integral over $T$ of the invariant $3$-form $\sum_j(-1)^jx_j\,dx_0\wedge\dots\wedge\widehat{dx_j}\wedge\dots\wedge dx_3/(x_0^2+\dots+x_3^2)^2$; in the chart $x_0 = 1$ this is $dz_1\wedge dz_2\wedge dz_3/(1+|z|^2)^2$, the volume form of $S^3$ in central projection. We shall not use this integral; the role it would play is taken by Schl\"afli's formula.
\end{rem}

\subsection{The Dehn element of a tetrahedron}\label{ssec:dehn}

Let $(V,Q,\omega)$ be an oriented quadratic space of rank $4$ over $E_m$, where $E$ contains $i$, and let $\Delta = ([v_1],\dots,[v_4])$ be a simplex in $\mathbb P(V)$, i.e.\ a tetrahedron. For each edge $\ell$ the plane $W_\ell$ is a quadratic space of rank $2$; we assume that it is \emph{orientable}, and we choose an orientation $\omega_\ell$. By \S\ref{ssec:quad} this means that the constant term of the Gram determinant of a basis of $W_\ell$ is a square in $E$. It is therefore a condition on $E$ as much as on $\Delta$: the same tetrahedron may fail to be orientable along its edges over $E$ and become so over a larger field, and it always does over a finite extension, since one only has to adjoin the six square roots in question. The Dehn element will turn out not to depend on the chosen $\omega_\ell$, so these orientations are not part of the data; what the field has to provide is their existence. This is how $E$ will be enlarged in \S\ref{ssec:descent}. We give $V/W_\ell$ the induced orientation $\omega^\ell$ of \eqref{eq:induced}. By Lemma~\ref{lem:link}, the edge $\Delta_\ell$ and its link $\Delta^\ell$ are simplices in the oriented planes $W_\ell$ and $V/W_\ell$, so their angles are defined.

\begin{definition}[Dehn element]\label{def:dehn}
The \emph{Dehn element} of $\Delta$ is
\[
  \Dehn(\Delta) := \sum_{\ell}\sgn(\ell)\;\mathcal Z_\ell\wedge\mathcal W_\ell\ \in\ \Lambda^2E_m^\times,\qquad \mathcal Z_\ell := \langle\Delta_\ell\rangle_{W_\ell,\omega_\ell},\quad \mathcal W_\ell := \langle\Delta^\ell\rangle_{V/W_\ell,\omega^\ell},
\]
the sum over the six edges $\ell\subseteq\{1,2,3,4\}$; thus $\mathcal Z_\ell = \langle v_a,v_b\rangle$ is the angle between the two vertices of the edge in the plane $W_\ell$, and $\mathcal W_\ell = \langle\bar v_c,\bar v_d\rangle$ the angle between the images of the other two vertices in the plane $V/W_\ell$, for $\ell = \{a<b\}$ and $\{1,2,3,4\}\setminus\ell = \{c<d\}$.
\end{definition}

As announced, the Dehn element does not depend on the orientations $\omega_\ell$: replacing $\omega_\ell$ by $-\omega_\ell$ replaces $\omega^\ell$ by $-\omega^\ell$, so it inverts both $\mathcal Z_\ell$ and $\mathcal W_\ell$, and $\mathcal Z_\ell^{-1}\wedge\mathcal W_\ell^{-1} = \mathcal Z_\ell\wedge\mathcal W_\ell$. It only depends on the points $[v_p]$. This is the definition in the form suggested by Goncharov's description of the Dehn invariant of a geodesic simplex as the sum, over its faces, of the face tensored with the spherical simplex cut out in the orthogonal space \cite[\S1.4, (3)]{G1}: in our case only the edges contribute, and the angle plays the role of both factors. The same definitions of faces and links make sense for simplices of any dimension.

\emph{In coordinates.} Let $V = E_m^4$ with the form $\pr{\ }{\ }_A$ and the orientation $\omega\det$, and write $\Dehn^{A,\omega}_m(v) := \Dehn(\Delta)$ for $\Delta = ([v_1],\dots,[v_4])$. We call $v = (v_1,\dots,v_4)$ a \emph{configuration} over $E_m$ for $A$ if $\Delta$ is a tetrahedron whose edge planes are orientable. For $\ell = \{a<b\}$ and $\{c<d\}$ its complement, put $c_\ell := \pr{v_a}{v_b}_A$, $\gamma_\ell := \pr{v_a}{v_a}_A\pr{v_b}{v_b}_A-c_\ell^2$, $s_\ell := \omega_\ell(v_a\wedge v_b)$, so that $s_\ell^2 = \gamma_\ell$, let $\beta_\ell$ be the minor of the Gram matrix with rows $a,b,c$ and columns $a,b,d$, and $\kappa_\ell := \omega\det(v_a,v_b,v_c,v_d)$. By Lemma~\ref{lem:link}, \eqref{eq:induced} and Lemma~\ref{lem:angle}(a), $Q(\bar v_c,\bar v_d) = \beta_\ell/\gamma_\ell$ and $\omega^\ell(\bar v_c\wedge\bar v_d) = \kappa_\ell/s_\ell$, so that
\begin{equation}\label{eq:ZW}
  \mathcal Z_\ell = \frac{c_\ell+is_\ell}{c_\ell-is_\ell},\qquad \mathcal W_\ell = \frac{\beta_\ell+is_\ell\kappa_\ell}{\beta_\ell-is_\ell\kappa_\ell}.
\end{equation}
In particular $\Dehn^{A,\omega}_m(v)$ only depends on the Gram matrix of $v_1,\dots,v_4$ and on $\omega\det(v_1,\dots,v_4)$.

For a real spherical tetrahedron, i.e.\ $E_m$ replaced by $\C$, $A = 1$, $\omega = 1$, and real vectors $v_p$ with $\det(v_1,\dots,v_4)>0$, Remark~\ref{rem:euclid} gives $\sgn(\ell)\,\mathcal Z_\ell\wedge\mathcal W_\ell = \exp(2i\theta_\ell)\wedge\exp(2i\psi_\ell)$, where $\theta_\ell$ is the spherical length of the edge $\ell$ and $\psi_\ell$ the dihedral angle along it. (Take $s_\ell>0$; then $\mathcal Z_\ell = \exp(2i\theta_\ell)$, and $\omega^\ell(\bar v_c\wedge\bar v_d)$ has the sign $\sgn(\ell)$, so $\mathcal W_\ell = \exp(2i\,\sgn(\ell)\psi_\ell)$.) So $\Dehn(\Delta)$ is the counterpart of $-4$ times the Dehn invariant $\sum_\ell\theta_\ell\otimes\psi_\ell$.

\begin{proposition}\label{prop:dehn}
Let $\Delta = ([v_1],\dots,[v_4])$ be a tetrahedron in $\mathbb P(V)$ as above.
\begin{enumerate}[nosep]
  \item[(a)] For $\rho\in S_4$, let $\Delta_\rho$ be the tetrahedron with the vertices permuted by $\rho$. Then $\Dehn(\Delta_\rho) = \sgn(\rho)\,\Dehn(\Delta)$.
  \item[(b)] Replacing $\omega$ by $-\omega$ changes the sign of $\Dehn(\Delta)$. If $g:V\to V'$ is an isometry onto another oriented quadratic space $(V',Q',\omega')$, then $\Dehn(g\Delta) = \pm\Dehn(\Delta)$, with the sign $+$ if and only if $g$ carries $\omega'$ to $\omega$.
  \item[(c)] If $\Delta$ is degenerate, then $\Dehn(\Delta) = 0$.
  \item[(d)] Let $[v_1],\dots,[v_5]\in\mathbb P(V)$ be such that for each $j$ the four points other than $[v_j]$, in increasing order, form a tetrahedron $\Delta_j$ as above. Then $\sum_{j=1}^5(-1)^j\,\Dehn(\Delta_j) = 0$.
\end{enumerate}
\end{proposition}

\begin{proof}
For any ordering $(a,b,c,d)$ of $\{1,2,3,4\}$, the element $\epsilon(a,b,c,d)\,\langle v_a,v_b\rangle\wedge\langle\bar v_c,\bar v_d\rangle$, where $\epsilon$ is the sign of the ordering and the angles are taken with respect to $\omega_{\{a,b\}}$ and the induced orientation, is unchanged under the transpositions $(ab)$ and $(cd)$ of the ordering, by Lemma~\ref{lem:angle}(b). Hence
\begin{equation}\label{eq:collapse}
  \Dehn(\Delta) = \frac14\sum_{(a,b,c,d)}\epsilon(a,b,c,d)\,\langle v_a,v_b\rangle\wedge\langle\bar v_c,\bar v_d\rangle,
\end{equation}
the sum over all orderings.

(a) Permuting the vertices by $\rho$ permutes the orderings in \eqref{eq:collapse}, multiplying their signs by $\sgn(\rho)$.

(b) Replacing $\omega$ by $-\omega$ replaces every $\omega^\ell$ by its negative and so inverts every $\mathcal W_\ell$. An isometry carrying $\omega'$ to $\omega$ transports the edge planes, their orientations and the induced orientations, so it preserves the Dehn element by Lemma~\ref{lem:angle}(d).

(c) Then $v_1\wedge\dots\wedge v_4 = 0$, so $\omega(v_1\wedge\dots\wedge v_4) = 0$, so $\omega^\ell(\bar v_c\wedge\bar v_d) = 0$ by \eqref{eq:induced}, and $\mathcal W_\ell = 1$ for every $\ell$.

(d) Choose orientations of the ten planes spanned by pairs of the $v_j$. Fix a pair $\{a,b\}$ and let $e<f<g$ be the remaining indices in $\{1,\dots,5\}$; the edge $\{a,b\}$ lies in $\Delta_e$, $\Delta_f$ and $\Delta_g$. In all three, $\langle v_a,v_b\rangle$ is the same, and the link is computed in the same oriented plane $V/W_{\{a,b\}}$. Let $s$ be the sign of the ordering $(a,b,e,f,g)$ of $\{1,\dots,5\}$. Moving the last entry to its place in the increasing order shows that
\[
  (-1)^g\epsilon_{\Delta_g}(a,b,e,f) = (-1)^f\epsilon_{\Delta_f}(a,b,g,e) = (-1)^e\epsilon_{\Delta_e}(a,b,f,g) = -s,
\]
where $\epsilon_{\Delta_j}$ is the sign of an ordering of the vertices of $\Delta_j$ relative to their increasing order; note that $(a,b,g,e,f)$ and $(a,b,f,g,e)$ also have sign $s$. By \eqref{eq:collapse}, the contribution of $\{a,b\}$ to $\sum_j(-1)^j\Dehn(\Delta_j)$ is therefore
\[
  -s\,\langle v_a,v_b\rangle\wedge\big(\langle\bar v_e,\bar v_f\rangle\,\langle\bar v_g,\bar v_e\rangle\,\langle\bar v_f,\bar v_g\rangle\big),
\]
and the product of the three angles is $1$ by Lemma~\ref{lem:angle}(b) and (c).
\end{proof}

\subsection{The Schl\"afli identity}\label{ssec:conf}

We first write down the ring over which the Dehn element is a universal object. Let
\[
  a = (a_{jj'})_{1\le j\le j'\le4},\qquad x = (x_{p,j})_{1\le p,j\le4}
\]
be ten and sixteen indeterminates, respectively, let $A$ be the symmetric matrix with entries $a_{jj'}$, and let $v_p := (x_{p,1},\dots,x_{p,4})$ for $1\le p\le4$. A point of ${\rm Spec}(E[a,x])$ with values in a ring $R\supseteq E$ is a pair $(A,v)$ consisting of a symmetric matrix $A\in M_4(R)$ and four vectors $v_1,\dots,v_4$ in $R^4$.

We want to invert exactly what is needed for $(A,v)$ to define a tetrahedron. By \S\ref{ssec:quad}, this means that $\det A$ is a unit, so that $\pr{\ }{\ }_A$ is non-degenerate, and that $\Gram\big((v_p)_{p\in F}\big)$ is a unit for every nonempty proper subset $F\subsetneq\{1,2,3,4\}$, so that the $v_p$, $p\in F$, span a non-degenerate subspace. Let $B_0$ be the localization of $E[a,x]$ at $\det A$ and at these fourteen Gram determinants, taken with respect to $\pr{\ }{\ }_A$. Thus, for a local ring $R\supseteq E$, a point of ${\rm Spec}(B_0)$ with values in $R$ is a pair $(A,v)$ which gives a non-degenerate quadratic space $\big(R^4,\pr{\ }{\ }_A\big)$ together with a tetrahedron $([v_1],\dots,[v_4])$ in $\mathbb P(R^4)$. In particular the $\gamma_\ell = \Gram(v_a,v_b)$ of \S\ref{ssec:dehn}, for $\ell = \{a<b\}$, are inverted in $B_0$.

It remains to add the square roots which the orientations require. Let
\[
  B := B_0[w,s_\ell]\big/\big(w^2-\det A,\ s_\ell^2-\gamma_\ell\big),
\]
with one variable $s_\ell$ for each of the six edges. A point of ${\rm Spec}(B)$ with values in $E_m$ gives in addition an orientation $\omega = w$ of $A$ and orientations $s_\ell$ of the six edge planes; its $v$ is then a configuration over $E_m$ for $A$ in the sense of \S\ref{ssec:dehn}. The elements $\mathcal Z_\ell$, $\mathcal W_\ell$ of \eqref{eq:ZW} are units of $B$. The \emph{Schl\"afli form} is
\[
  \Omega := \sum_\ell\sgn(\ell)\,d\log\mathcal Z_\ell\wedge d\log\mathcal W_\ell\ \in\ \Omega^2_{B/E}.
\]
For a real spherical tetrahedron, $\sgn(\ell)\,d\log\mathcal Z_\ell\wedge d\log\mathcal W_\ell = (2i\,d\theta_\ell)\wedge(2i\,d\psi_\ell)$, so that
\[
  \Omega = -4\sum_\ell d\theta_\ell\wedge d\psi_\ell = -4\,d\Big(\sum_\ell\theta_\ell\,d\psi_\ell\Big) = -8\,d(dV) = 0
\]
by Schl\"afli's formula $\sum_\ell\theta_\ell\,d\psi_\ell = 2\,dV$ for the volume $V$ of a spherical tetrahedron \cite{Mi}; see \cite[\S4.8]{G1} for its Hodge-theoretic interpretation. The following proposition is the algebraic form of this computation.

\begin{proposition}[Schl\"afli identity]\label{prop:schlafli}
$\Omega = 0$ in $\Omega^2_{B/E}$.
\end{proposition}

\begin{proof}
The ring $B$ is obtained from $B_0$ by adjoining square roots of units. So it is finite \'etale over $B_0$, free with basis the monomials $w^{e_0}\prod_\ell s_\ell^{e_\ell}$, $e\in\{0,1\}^7$, and $\Omega^2_{B/E} = B\otimes_{B_0}\Omega^2_{B_0/E}$. The group $G = \{\pm1\}^7$ acts on $B$ by changing the signs of $w$ and the $s_\ell$, and $B^G = B_0$, since $2$ is invertible. Changing the sign of $s_\ell$ inverts $\mathcal Z_\ell$ and $\mathcal W_\ell$, and changing the sign of $w$ inverts all the $\mathcal W_\ell$ (Lemma~\ref{lem:angle}(b)). Hence $w\,\Omega$ is invariant under $G$, so it lies in $\Omega^2_{B_0/E}$. Its coefficients with respect to the basis of $\Omega^2_{B_0/E}$ given by the $da_{jj'}\wedge dx_{p,j}$, etc., are rational functions in $(a,x)$ with coefficients in $\Q(i)$. It therefore suffices to show that they vanish on a nonempty open subset of the real points $(a,x)\in\mathbb R^{10}\times\mathbb R^{16}$.

Let $U$ be the open set where $A$ is positive definite and $\det(v_1,\dots,v_4)>0$. On $U$ the Gram determinants of all subsets of the vertices are positive, because the Gram matrix of linearly independent vectors in a positive definite space is positive definite; in particular $U$ lies in the real points of ${\rm Spec}(B_0)$, and $\gamma_\ell>0$. At a point of $U$ choose $w = \sqrt{\det A}>0$ and $s_\ell = \sqrt{\gamma_\ell}>0$. Then $v_1,\dots,v_4$ are the vertices of a spherical tetrahedron in the unit sphere of the Euclidean space $(\mathbb R^4,A)$, and as above $\mathcal Z_\ell = \exp(2i\theta_\ell)$ and $\mathcal W_\ell = \exp(2i\,\sgn(\ell)\psi_\ell)$, where $\theta_\ell$, $\psi_\ell$ and the volume $V$ are smooth functions on $U$. Schl\"afli's formula $2\,dV = \sum_\ell\theta_\ell\,d\psi_\ell$ holds on $U$: for fixed $A$ it is the classical formula, and in general all three of $V$, $\theta_\ell$ and $\psi_\ell$ are pulled back from the case $A = 1$ by the map $(A,v)\mapsto A^{1/2}v$. Hence $\Omega = -8\,d(dV) = 0$ at the points of $U$, so the coefficients of $w\,\Omega$ vanish on $U$ and therefore $w\,\Omega = 0$ in $\Omega^2_{B_0/E}$. Since $w$ is a unit, $\Omega = 0$ in $\Omega^2_{B/E}$.
\end{proof}

\subsection{\texorpdfstring{The infinitesimal volume over $E$}{The infinitesimal volume over E}}\label{ssec:vol}

Let $n\ge2$, let $(V,Q,\omega)$ be an oriented quadratic space of rank $4$ over $E_n$, and let $\Delta$ be a tetrahedron in $\mathbb P(V)$ whose edge planes are orientable. A \emph{lift} of $(V,Q,\omega,\Delta)$ to $E_{2n}$ is an oriented quadratic space $(\tilde V,\tilde Q,\tilde\omega)$ over $E_{2n}$ together with a tetrahedron $\tilde\Delta$ in $\mathbb P(\tilde V)$ and an isomorphism $\tilde V\otimes_{E_{2n}}E_n\simeq V$ which carries $\tilde Q$, $\tilde\omega$ and $\tilde\Delta$ to $Q$, $\omega$ and $\Delta$.

Lifts exist, and in coordinates they are easy to describe. Choose a basis of $V$, so that $V = E_n^4$, $Q = \pr{\ }{\ }_A$, the orientation is $\omega\det$ with $\omega\in E_n$, and $\Delta = ([\overline v_1],\dots,[\overline v_4])$ for a configuration $\overline v$ over $E_n$ for $A$ (\S\ref{ssec:dehn}). Let $\tilde A\in M_4(E_{2n})$ be a symmetric lift of $A$, and $v_p\in E_{2n}^4$ lifts of the $\overline v_p$. There is a unique lift $\tilde\omega$ of $\omega$ with $\tilde\omega^2 = \det\tilde A$, and $\big(E_{2n}^4,\pr{\ }{\ }_{\tilde A},\tilde\omega\det,([v_1],\dots,[v_4])\big)$ is a lift: being a tetrahedron, and having orientable edge planes, are conditions on constant terms only. We call $(\tilde A,\tilde\omega,v)$ a lift of $(A,\omega,\overline v)$. Conversely, every lift is of this form, since any lift of the basis of $V$ to $\tilde V$ is again a basis. The orientations of the edge planes lift uniquely as well.

\begin{definition}\label{def:Vol}
The \emph{infinitesimal volume} of $\Delta$ is
\[
  \Vol^\omega_n(\Delta) := -\frac14\,d^{-1}\circ\La_n\big(\Dehn(\tilde\Delta)\big)\ \in\ D_1(E_n)
\]
for a lift $(\tilde V,\tilde Q,\tilde\omega,\tilde\Delta)$ of $(V,Q,\omega,\Delta)$. In coordinates as above we write $\Vol^{A,\omega}_n(\overline v) := \Vol^\omega_n(\Delta)$, so that $\Vol^{A,\omega}_n(\overline v) = -\frac14\,d^{-1}\circ\La_n\big(\Dehn^{\tilde A,\tilde\omega}_{2n}(v)\big)$ for a lift $(\tilde A,\tilde\omega,v)$ of $(A,\omega,\overline v)$; for $A = 1$ and $\omega = 1$ we write $\Vol_n(\overline v)$.
\end{definition}

The factor $-\frac14$ compensates for the factor $-4$ relating the Dehn element to $\sum_\ell\theta_\ell\otimes\psi_\ell$; see Remark~\ref{rem:jets}.

\begin{theorem}\label{thm:vol}
$\Vol^\omega_n(\Delta)$ does not depend on the choice of the lift.
\end{theorem}

\begin{proof}
Choose a basis of $V$. By the above, two lifts are given in coordinates by two lifts $(\tilde A,\tilde\omega,v)$ and $(\tilde A',\tilde\omega',v')$ of $(A,\omega,\overline v)$, and by Proposition~\ref{prop:dehn}(b) their Dehn elements are $\Dehn^{\tilde A,\tilde\omega}_{2n}(v)$ and $\Dehn^{\tilde A',\tilde\omega'}_{2n}(v')$. Together with the lifted orientations of the edge planes, these are $E_{2n}$-valued points of ${\rm Spec}(B)$ with the same reduction. So Lemma~\ref{lem:var}, over the base field $E$, applies with $D$ the difference of the corresponding maps $B\to E_{2n}$. Taking $f := \mathcal Z_\ell$, $g := \mathcal W_\ell$ and summing over $\ell$, we see that the two values of $\La_n(\Dehn_{2n})$ differ by $\overline\varphi_*\big(\iota_D(\Omega)\big) = 0$, by Proposition~\ref{prop:schlafli}.
\end{proof}

The proof is parallel to the argument in Remark~\ref{rem:dilog}. There $B = k[z,\frac{1}{z(1-z)}]$, the element is $\delta[z] = (1-z)\wedge z$, and the $2$-form $d\log(1-z)\wedge d\log z$ vanishes because $\Omega^2_{B/k} = 0$. Here the element is the Dehn element, and the vanishing of the $2$-form is the Schl\"afli identity.

\begin{theorem}\label{thm:scissors}
Let $\Delta$ be a tetrahedron in $\mathbb P(V)$ whose edge planes are orientable.
\begin{enumerate}[nosep]
  \item[(R1)] $\Vol^\omega_n(\Delta_\rho) = \sgn(\rho)\,\Vol^\omega_n(\Delta)$ for $\rho\in S_4$.
  \item[(R2)] If $g:V\to V'$ is an isometry onto an oriented quadratic space $(V',Q',\omega')$ carrying $\omega'$ to $\omega$, then $\Vol^{\omega'}_n(g\Delta) = \Vol^\omega_n(\Delta)$. Moreover $\Vol^{-\omega}_n = -\Vol^\omega_n$. In particular $\Vol^\omega_n(g\Delta) = \det(g)\,\Vol^\omega_n(\Delta)$ for $g\in O(V)$.
  \item[(R3)] $\Vol^\omega_n(\Delta) = 0$ if $\Delta$ is degenerate.
  \item[(R4)] For points $x_1,\dots,x_5$ of $\mathbb P(V)$ as in Proposition~\ref{prop:dehn}(d), with tetrahedra $\Delta_j$ whose edge planes are orientable, $\sum_{j=1}^5(-1)^j\,\Vol^\omega_n(\Delta_j) = 0$.
\end{enumerate}
\end{theorem}

The relations (R2), (R3) and (R4) are the analogues of the relations (3.1)(iii), (i) and (ii) by which Dupont and Sah present the scissors congruence group of totally asymptotic hyperbolic polytopes in terms of tuples of vertices \cite[\S3]{DS}: invariance under isometries $g$ up to the sign $\det(g)$, vanishing on degenerate tuples, and, in dimension three, the alternating sum over five points. The same relations, for simplices described by their faces, define Goncharov's group $S_2(\C)$ \cite[\S3.4, R0)--R3)]{G1}. Our tetrahedra are generic: all their vertices, edges and faces are non-degenerate. The classical cutting relation $[x_0,x_1,x_2,x_3] = [x_0,x_4,x_2,x_3]+[x_4,x_1,x_2,x_3]$, for $x_4$ on the line through $x_0$ and $x_1$, is formally the case of (R4) in which two of the five terms are degenerate. It is not available here: the triple $x_0,x_1,x_4$ only spans a plane, so the two quadruples containing it are not tetrahedra in our sense.

\begin{proof}
By Theorem~\ref{thm:vol} we may compute with any lift, and $d^{-1}\circ\La_n$ is $\Q$-linear.

(R1) Apply Proposition~\ref{prop:dehn}(a) to a lift.

(R2) Let $(\tilde V,\tilde Q,\tilde\omega,\tilde\Delta)$ be a lift of $(V,Q,\omega,\Delta)$. Composing the identification $\tilde V\otimes E_n\simeq V$ with $g$ makes it a lift of $(V',Q',\omega',g\Delta)$, with the same Dehn element. The second statement follows from Proposition~\ref{prop:dehn}(b). For $g\in O(V)$, we have $\det(g) = \pm1$, and $g$ carries $\det(g)\,\omega$ to $\omega$; so $\Vol^\omega_n(g\Delta) = \det(g)\Vol^{\det(g)\omega}_n(g\Delta) = \det(g)\Vol^\omega_n(\Delta)$.

(R3) By (R1) we may assume $\Delta = ([\overline v_1],\dots,[\overline v_4])$ with $\overline v_4 = \sum_{j\le3}\lambda_j\overline v_j$, $\lambda_j\in E_n$. Choose a lift $\tilde V$ of $V$, lift the $\lambda_j$ and $\overline v_1,\overline v_2,\overline v_3$ arbitrarily, and put $v_4 := \sum_j\tilde\lambda_jv_j$. This gives a degenerate lift of $\Delta$; apply Proposition~\ref{prop:dehn}(c).

(R4) Lift the five points to one lift $\tilde V$ of $V$; the four-element subsets then give lifts of the $\Delta_j$. Apply Proposition~\ref{prop:dehn}(d).
\end{proof}

\subsection{\texorpdfstring{Descent to $k$}{Descent to k}}\label{ssec:descent}

Fix an algebraic closure $\bar k$ of $k$; all extensions of $k$ below are taken inside $\bar k$. Let $(V,Q)$ be a quadratic space of rank $4$ over $k_n$. The constant term of the Gram determinant of a basis of $V$ lies in $k^\times$, and a change of basis multiplies it by a square. Let $\disc(V)\in k^\times/k^{\times2}$ be its class, the \emph{discriminant} of $V$, let $\sqrt{\disc(V)}\in\bar k$ be a square root of a representative, and $k' := k\big(\sqrt{\disc(V)}\big)$, which only depends on $V$. Then $V\otimes_{k_n}k'_n$ has an orientation $\omega$, unique up to sign, and $\omega$ is defined over $k_n$ if and only if $\disc(V)$ is trivial, i.e.\ $k' = k$.

\begin{definition}\label{def:conf-k}
Let $\Delta$ be a tetrahedron in $\mathbb P(V)$. Let $E\subseteq\bar k$ be a finite extension of $k'$ containing $i$ and square roots of the constant terms of the Gram determinants of the edges of $\Delta$. We define $\Vol^\omega_n(\Delta)\in D_1(E_n)$ by Definition~\ref{def:Vol}, applied to $(V,Q,\omega,\Delta)\otimes_{k_n}E_n$. In coordinates, i.e.\ for $V = k_n^4$ with $Q = \pr{\ }{\ }_A$, the orientation $\omega\det$, and $\Delta = ([\overline v_1],\dots,[\overline v_4])$, we call $\overline v$ a \emph{configuration over $k_n$ for $A$} and write $\Vol^{A,\omega}_n(\overline v) := \Vol^\omega_n(\Delta)$.
\end{definition}

\begin{theorem}\label{thm:descent}
Let $\Delta$ be a tetrahedron in $\mathbb P(V)$.
\begin{enumerate}[nosep]
  \item[(i)] $\Vol^\omega_n(\Delta)$ does not depend on the choice of $E$, and lies in $\sqrt{\disc(V)}\cdot D_1(k_n)$. If $\disc(V)$ is non-trivial, the non-trivial automorphism of $k'/k$ maps $\Vol^\omega_n(\Delta)$ to $-\Vol^\omega_n(\Delta)$.
  \item[(ii)] The relations (R1)--(R4) of Theorem~\ref{thm:scissors} hold for tetrahedra in $\mathbb P(V)$, with isometries over $k_n$ in (R2), tetrahedra degenerate over $k_n$ in (R3), and points $x_1,\dots,x_5$ of $\mathbb P(V)$ over $k_n$ in (R4), any four of which form a tetrahedron.
\end{enumerate}
Consequently $\Vol^\omega_n$ induces a homomorphism
\[
  \Vol^\omega_n:\ \mathcal P(V)\longrightarrow D_1(k'_n),
\]
where $\mathcal P(V)$ is the abelian group generated by the symbols $[\Delta]$, for the tetrahedra $\Delta$ in $\mathbb P(V)$, subject to the relations
\[
  [\Delta_\rho] = \sgn(\rho)\,[\Delta],\quad [g\Delta] = \det(g)\,[\Delta],\quad [\Delta] = 0\ \text{if $\Delta$ is degenerate},\quad \sum_{j=1}^5(-1)^j[\Delta_j] = 0,
\]
for $\rho\in S_4$, $g\in O(V)$, and points $x_1,\dots,x_5$ of $\mathbb P(V)$ as in (ii). It is a scissors congruence group of generic tetrahedra. For $V = k_n^4$ with the standard form and $\omega = \det$ we write $\mathcal P(k_n)$ and obtain $\Vol_n:\mathcal P(k_n)\to D_1(k_n)$.
\end{theorem}

\begin{proof}
(i) Choose a basis of $V$, so that $V = k_n^4$, $Q = \pr{\ }{\ }_A$ with $A\in M_4(k_n)$, the orientation is $\omega\det$ with $\omega\in k'_n$, $\omega^2 = \det A$, and $\Delta$ is given by a configuration $\overline v$ over $k_n$ for $A$. If $E\subseteq E'$, the construction is compatible with the inclusions $E_m\subseteq E'_m$ and $D_1(E_n)\subseteq D_1(E'_n)$. Since any two admissible fields lie in a common one, $\Vol^{A,\omega}_n(\overline v)$ does not depend on $E$, and we may take $E$ to be Galois over $k$. Let $\tau\in\operatorname{Gal}(E/k)$, acting on $E_m$ through the coefficients. It commutes with $\log^\circ$ and with $d:E_m\to\Omega^1_{E_m/E} = E_m\,dt/(t^{m-1}dt)$, on which it also acts through the coefficients; hence $\tau\circ d^{-1}\circ\La_n = d^{-1}\circ\La_n\circ\Lambda^2\tau$. Choose a lift $(\tilde A,\tilde\omega,v)$ with $\tilde A$ and $v$ defined over $k_{2n}$, so that $\tilde\omega\in k'_{2n}$. Then $\tau$ fixes $\tilde A$ and $v$, maps the orientations of the edge planes to orientations, maps $i$ to $\pm i$, and maps $\tilde\omega$ to $\tilde\omega$ or $-\tilde\omega$ according as $\tau$ is trivial on $k'$ or not. By Lemma~\ref{lem:angle}(b), the Dehn element does not depend on the orientations of the edge planes and is unchanged when $i$ is replaced by $-i$, and by Proposition~\ref{prop:dehn}(b) it changes sign with $\tilde\omega$. Hence $\tau\big(\Vol^{A,\omega}_n(\overline v)\big) = \pm\Vol^{A,\omega}_n(\overline v)$, with the sign $+$ if and only if $\tau$ is trivial on $k'$. Since $D_1(E_n) = \oplus_wEt^w$ with $\operatorname{Gal}(E/k)$ acting on the coefficients, taking invariants under $\operatorname{Gal}(E/k')$ gives $\Vol^{A,\omega}_n(\overline v)\in D_1(k'_n)$. If $k'\ne k$, then $D_1(k'_n) = D_1(k_n)\oplus\sqrt{\disc(V)}\,D_1(k_n)$, and the non-trivial automorphism of $k'/k$ acts by $1$ and $-1$ on the two summands.

(ii) Choose $E$ such that the edge planes of all the tetrahedra involved are orientable over $E_n$. Then (ii) follows from Theorem~\ref{thm:scissors} and (i).
\end{proof}

\begin{rem}[The case $n = 2$]\label{rem:n2}
Here $D_1(k_2) = kt^3$. Write $\log^\circ(a) = a_1t+a_2t^2+\cdots$ and $\log^\circ(b) = b_1t+b_2t^2+\cdots$. Since $t\,dt = 0$ in $\Omega^1_{E_2/E}$, we have $\La_2(a\wedge b) = (a_1b_2-a_2b_1)\,dt\otimes t^2$, and $d^{-1}(t^2dt) = t^3/3$. Hence
\[
  \Vol^{A,\omega}_2(\overline v) = -\frac{t^3}{12}\sum_\ell\sgn(\ell)\,(a_{\ell,1}b_{\ell,2}-a_{\ell,2}b_{\ell,1}),
\]
where $\log^\circ\mathcal Z_\ell(v) = \sum_ja_{\ell,j}t^j$ and $\log^\circ\mathcal W_\ell(v) = \sum_jb_{\ell,j}t^j$, for any lift $v$ to $E_3$; this suffices, since only coefficients of $t$-degree $\le2$ enter. By Theorem~\ref{thm:vol}, this does not depend on the $t^2$-coefficients of the lift. It is a function of the $1$-jet $\overline v$, cubic in its $t$-coefficients, with coefficients rational functions of the constant terms. For example, for $A = 1$ and $\overline v_p = e_p+b_pt$, with $e_p$ the standard basis, $b_1 = (0,1,1,0)$, $b_2 = (0,0,1,1)$ and $b_3 = b_4 = 0$, one finds $\Vol_2(\overline v) = -t^3$.
\end{rem}

\begin{rem}[Angle jets]\label{rem:jets}
Let $v$ be a lift of $\overline v$ to $E_{2n}$. Write formally $\log^\circ\mathcal Z_\ell(v) = 2i\theta_\ell$ and $\sgn(\ell)\log^\circ\mathcal W_\ell(v) = 2i\psi_\ell$, with $\theta_\ell,\psi_\ell\in tE_{2n}$ the jets of the edge lengths and the dihedral angles along $v$, and let $d\theta_\ell$, $d\psi_\ell\in\Omega^1_{E_n/E}$ be the differentials of their reductions to $E_n$. Then
\[
  \Vol^{A,\omega}_n(\overline v) = d^{-1}\Big(\sum_\ell d\theta_\ell\otimes(\psi_\ell)_{\ge n}-d\psi_\ell\otimes(\theta_\ell)_{\ge n}\Big),
\]
an integrated, truncated version of the Schl\"afli form $\sum_\ell\theta_\ell\,d\psi_\ell$. Lemma~\ref{lem:explicit} expresses $\Vol_n^{A,\omega}(\overline v)$ in terms of the coefficients of these jets. The construction is natural for the automorphisms $\sigma_\lambda:t\mapsto\lambda t$, $\lambda\in k^\times$: the automorphism $\sigma_\lambda$ of $E_{2n}$ carries a lift of $(A,\omega,\overline v)$ to a lift of its base change along $\sigma_\lambda$, and it commutes with $\log_{\ge n}$ and with $d$, hence with $\La_n$ and $d^{-1}$. So $\Vol^{\sigma_\lambda A,\sigma_\lambda\omega}_n(\sigma_\lambda\overline v) = \sigma_\lambda\big(\Vol^{A,\omega}_n(\overline v)\big)$; for $A = 1$ the coefficient of $t^w$ is homogeneous of degree $w$ in the jets of $\overline v$. This proves part (iii) of Theorem~A. For example, for $A = 1$, $n = 3$ and $\overline v_p = e_p+b_pt+c_pt^2$ with $b_1 = (0,1,1,0)$, $b_2 = (0,0,1,1)$, $b_3 = (1,0,0,1)$, $b_4 = (1,1,0,0)$ and $(c_1,c_2,c_3,c_4) = (e_1,e_4,e_2,e_3)$, one finds $\Vol_3(\overline v) = 5t^4+4t^5$.
\end{rem}

\begin{rem}[Hyperbolic tetrahedra]\label{rem:hyp}
For $A = \operatorname{diag}(1,1,1,-1)$ we have $\disc(V) = -1$, $k' = k(i)$ and $\omega = \pm i$. Configurations for $A$ are the algebraic counterparts of tetrahedra in hyperbolic $3$-space in the projective model, the vertices lying inside the quadric $x_1^2+x_2^2+x_3^2 = x_4^2$, and Theorem~\ref{thm:descent} gives their infinitesimal volume, with values in $i\cdot D_1(k_n)$. Indeed, for the vertices of a real hyperbolic tetrahedron with $\det(v_1,\dots,v_4)>0$ and $\omega = i$, one has $\sgn(\ell)\,\mathcal Z_\ell\wedge\mathcal W_\ell = \exp(2\epsilon_\ell\theta_\ell)\wedge\exp(2i\epsilon_\ell\psi_\ell)$, where $\theta_\ell$ is the hyperbolic length of $\ell$, $\psi_\ell$ the dihedral angle along it, and $\epsilon_\ell = \pm1$ depends on the orientation of the edge plane. So the Dehn element is the counterpart of $4i$ times the hyperbolic Dehn invariant, just as it is the counterpart of $-4$ times the spherical one. Algebraically the two cases differ only by the choice of the form; compare Goncharov's construction, in which spherical and hyperbolic simplices both give elements of the same group $S_2(\C)$ \cite[\S\S1.5, 3.4]{G1}. The choice of $\omega$ corresponds to the choice of a family of maximal isotropic subspaces of the quadric in \cite[\S3.4, R2)b)]{G1}. Ideal vertices, with $\pr{v_p}{v_p}_A = 0$, are excluded, since the vertices of a tetrahedron are non-isotropic.
\end{rem}

\begin{rem}[Relation to the Bloch group]\label{rem:bloch}
The Bloch group gives an analogue of the volume map for pairs of triangles over $k_n$: map a pair of triangles to $B_2(k_n)$ as in \cite{BGSV} and apply the additive dilogarithm \cite[\S1.3]{U1}. The construction here is the three-dimensional counterpart, taken from the Dehn side: $\Vol_n$ is $d^{-1}\circ\La_n$ applied to the Dehn element of a lift, just as $\Lg$ is $d^{-1}\circ\La_n$ applied to $\delta$ of a lift. One expects a natural map $\mathcal P(k_n)\to B_2(k_n)$ under which the Dehn element of a lift corresponds to $\delta$ of a lift, and $\Vol_n$ corresponds to a multiple of $\Lg$. This would be an analogue over $k_n$ of Goncharov's conjectural homomorphism $s_2:S_2(F)\to B_2(F)$ \cite[Conjecture 6.2]{G1}. For $\mathcal P(S^3)$, a map of this kind to the Bloch complex of $\C$ was constructed by Dupont, Parry and Sah using the Hopf fibration \cite{DPS}, \cite[\S6.1]{G1}. Theorem~\ref{thm:vol} does not use such a map; its proof uses only the Schl\"afli identity. In analogy with this and with Sydler's theorem, one may also expect that the class of a tetrahedron in $\mathcal P(k_n)$, modulo the classes of tetrahedra over $k$, is determined by its Dehn element $\Dehn^{1,1}_n(\overline v)\in\Lambda^2E_n^\times$, computed without lifting, together with $\Vol_n(\overline v)$.
\end{rem}

\medskip
\noindent\textbf{Acknowledgements.} The author used Claude (Anthropic) for assistance with the exposition and for numerical checks of some of the identities; the author takes full responsibility for the content of the paper.

\bigskip
\noindent\textsc{Ko\c{c} University, Mathematics Department, Rumelifeneri Yolu, 34450, Istanbul, Turkey}\\
\textit{E-mail address:} \texttt{sunver@ku.edu.tr}

\end{document}